  \documentclass[12pt,twoside]{article}
  \usepackage{amsfonts}

   \def\R{\mathbb{R}}
   \def\N{\mathbb{N}}
   
   \def\e{{\varepsilon}}
   \def\D{{\nabla}}
   \def\cA{{\cal A}}

   \def\cG{{\cal G}}

   \def\cJ{{\cal J}}    
   \def\cK{{\cal K}}

   \def\cS{{\cal S}}

   \def\supp{\mathop{\rm supp}\nolimits}
   
   \def\dist{\mathop{\rm dist}\nolimits}
   
   \def\const{\mathop{\rm const}\nolimits}
   
   \def\loc{\mathop{\rm loc}\nolimits}

\def\qed{\hfill{\em q.e.d.}}
\def\proof{\noindent{\underline {\sf Proof} \hspace{2mm}}}

\newtheorem{df}{Definition}[section]
\newtheorem{prop}[df]{Proposition}
\newtheorem{lemma}[df]{Lemma}
\newtheorem{teo}[df]{Theorem}
\newtheorem{rem}[df]{Remark}

 \newcommand{\sezione}[1]{\section{#1}\setcounter{equation}{0}}

\begin{document}


   \title{Semilinear elliptic problems in unbounded domains with
 unbounded boundary\thanks{Work supported by  
the national research  project ``Metodi 
variazionali e topologici nello studio di fenomeni non lineari".}
\author{
Riccardo Molle\thanks{Dipartimento di
Matematica, Universit\`a di Roma ``Tor Vergata", Via della Ricerca
Scientifica n$^o$ 1 - 00133 Roma. E-mail: {\sf  molle@mat.uniroma2.it}}}}
\date{}
  \maketitle


\begin{abstract}

This paper deals with a class of singularly perturbed nonlinear elliptic 
problems $(P_\e)$ with subcritical nonlinearity.	
	The coefficient of the linear part is assumed to concentrate 
in a point of the domain, as $\e\to 0$, and the domain is supposed to
be unbounded and with unbounded boundary.
	Domains that enlarge at infinity, and whose boundary flattens
or shrinks at infinity, are considered.
	It is proved that in such domains problem $(P_\e)$ has at
least 2 solutions.

\end{abstract}


{\small
\noindent{\bf Key words:}
Unbounded domains. Unbounded boundary. Concentrating
potential. Multiple solutions.


\noindent{\bf A.M.S. subject classification 2000:} 35J60, 35J20, 35J25.
}


\sezione{Introduction}


This paper deals with the problem
\[
(P_\e) \hspace{1cm}
\left\{
        \begin{array}{ll}
         -\e^2\Delta u+a_\e(x)u=u^{p-1} & \mbox{ in }\Omega\\
         u>0                            & \mbox{ in }\Omega\\
         u=0                            & \mbox{ on }\partial\Omega,
        \end{array}
\right.
\]
where $\Omega$ is an unbounded domain in $\R^N$, $N\geq 2$, having smooth
boundary, $\e>0$, $p>2$ and $p<2N/(N-2)$ when $N\geq 3$.
The function $a_\e(x)$ is assumed to concentrate in a point
$x_0\in\Omega$, when $\e$ goes to 0; that is $a$ is assumed to have
the following form
\begin{equation}
\label{a}
a_\e(x)=a_0+\alpha\left(\frac{x-x_0}{\e}\right)
\end{equation}
with
\begin{equation}
\label{a1}
x_0\in\Omega, \ a_0>0,\ \alpha\in L^{N/2}(\R^N), \
\alpha\geq 0, \ |\alpha|_{L^{N/2}(\R^N)}>0,
\end{equation}
\begin{equation}
\label{a2}
\int_{\R^N}\alpha (x)e^{2|x|}\left(1+|x|^{{N-1\over
2}\sigma}\right)dx<\infty \ \mbox{ for some }\ \sigma\in(1,2].
\end{equation}

\vspace{1mm}

There is a large literature on problem $(P_\e)$ in the case
$\Omega=\R^N$ (Schr\"o\-dinger equations). The existence and
multiplicity of solutions is related to some critical points of
$a_\e(x)$ or to some topological property of the sublevels of
$a_\e(x)$ (see, for example, \cite{ABCi,FW,MMP,R} and references
therein).

When $\Omega$ is an exterior domain,
i.e. $\Omega=\R^N\setminus\bar\omega$, with $\omega\subset\R^N$ bounded,
and $a_\e(x)\equiv\const$ we refer to \cite{BL} and \cite{BC}, where
the existence of at least one solution $u_\omega$ for problem $(P_\e)$
is proved. 
	In \cite{MP98} the behaviour of the ``energy'' of the
solution $u_\omega$ is studied, as the ``hole'' $\omega$
increase. Taking into account this behaviour, a multiplicity result is
found in \cite{MP00}, when $\omega$ has several suitable connected
components (see also \cite{GP} for a multiplicity result).
If $\Omega$ is an exterior domain and $a_\e(x)$ concentrates in a
point of $\Omega$ (see (\ref{a})-(\ref{a2})), then in
\cite{CM} it  has been proved that problem $(P_\e)$ has at least three
solutions.

In this paper we are interested in problem $(P_\e)$ when not only
$\Omega$, but also $\partial \Omega$ is unbounded.
	In the autonomous case this problem has been studied in 
\cite{EL} and a non-existence result has been proved for a class of
domains, that includes half-spaces.
	On the other hand, if $\Omega$ is a strip-like domain  in
\cite{E} it has been proved that problem $(P_\e)$ has a solution, 
found taking advantage of the symmetry properties of the domain.
	We mention also that an existence result has been
proved in \cite{G},  assuming  $\R^N\setminus\Omega$ of
``small capacity''. 

In the present work we consider problem $(P_\e)$ in unbounded domains 
 that do not enjoy of symmetry properties and whose
complement is not required to be small.
	Of course we have to impose some restriction to the shape of
$\Omega$.
	For $x\in\Omega$ and $y\in\partial\Omega$, we set
\[
r(x)=\sup\{\rho>0\ : \ \exists\bar x\in\Omega\ \mbox{ such that } x\in
B(\bar x,\rho)\mbox{ and } B(\bar x,\rho)\subset\Omega\},
\]
\begin{eqnarray*}
h(y)=\sup\{\dist(z,T_{\partial\Omega,y}\cap B(y,1))&:& z\mbox{ is in
the connected component}\\
& & \mbox{ of }\partial\Omega\cap B(y,1)\mbox{ containing } y\}
\end{eqnarray*}
where $B(y,r)$, $r>0$, denotes the ball centered in $y$ and with radius $r$ and
$T_{\partial\Omega,y}$ is the hyperplane tangent to $\partial\Omega$
in $y$; then we assume that $\Omega$ satisfies
\[
\begin{array}{ll}
\vspace{2mm}
{\displaystyle \lim_{R\to \infty} \inf\{r(x)\ :\ x\in\Omega,\
|x|=R\}=+\infty,}& \hspace{3cm} (C_1)\\
{\displaystyle \lim_{R\to \infty}\sup\{ h(y)\ :\
y\in\partial\Omega,\ |y|=R\}=0}. & \hspace{3cm} (C_2)
\end{array}
\]

Assumption $(C_1)$ implies that the domain $\Omega$ enlarges at
infinity while by assumption $(C_2)$ its boundary either flattens or shrinks,
at infinity.
	For example, domains that verify assumptions $(C_1)$ and
$(C_2)$ are
\[
\begin{array}{l}
\vspace{2mm}
\Omega_1^\alpha=\left\{(x',x_N)\in\R^N\ :\ x'\in\R^{N-1},\ x_N\in\R,\ 
|x'|^2>{1\over(1+x_N^2)^\alpha}\right\},\quad\alpha>0,\\
\vspace{2mm}
\Omega'=\{(x',x_N)\in\R^N\ :\ x'\in\R^{N-1},\ x_N\in\R,\ 
x_N>(1+|x'|^2)^{1/2}\},\\
\Omega{''}=\R^N\setminus\bar\Omega' 
\end{array}
\]
and $\Omega_2^\alpha=\R^N\setminus\bar\omega^\alpha$, where, for
$\alpha>{1\over 2}$, 
\[
\omega^\alpha=
\left\{(x',x_N)\in\R^N\ :\ x'\in\R^{N-1},\ x_N\in\R,\ |x'|>1,\ 
x_N<{(|x'|-1)^{1/2}\over|x'|^\alpha}\right\}.
\]
Let us remark that even if the sets $\Omega_1^\alpha$ have the complementary
``small'', they are more general than the domains considered in \cite{G},
where it is required that $\alpha>{1\over N-3}$, when $N>3$.

We are fronting problem $(P_\e)$ in a variational way. An essential
difficulty in this approach is caused by the lack of compactness
due to the non-compact embedding $H_0^1(\Omega)\hookrightarrow
L^P(\Omega)$, when $\Omega$ is unbounded.
	If $\Omega=\R^N$ or $\Omega$ is an exterior domain, then, by the
concentration-compactness principle (see \cite{BC}, \cite{L}), a local
compactness condition holds, that allows to apply mini-max techniques
in some energy intervals and permit to give existence and multiplicity results (see \cite{MMP},
\cite{CM}).
	When $\partial\Omega$ is unbounded the compactness situation
is worse, in general. Indeed, as we shall see in detail in Remark
\ref{SL}, when $\Omega$ is either a strip-like domain or the exterior
of a cylinder, the local compactness condition fails.
	Nevertheless we prove that, if the assumptions $(C_1)$ and
$(C_2)$ are satisfied, then the local compactness is restored in some
sense and we can  apply  topological methods
in order to obtain solutions for problem $(P_\e)$.
	The result we obtain is the following:

\begin{teo}
\label{T}
Assume that $a_\e$ is of the form (\ref{a}) and satisfies (\ref{a1}), (\ref{a2})
and suppose that $\Omega$ verifies assumptions $(C_1)$ and $(C_2)$; 
then there exists
$\bar\e>0$ such that for $\e\in(0,\bar\e)$ problem $(P_\e)$
has at least 2 solutions.
\end{teo}

One of the solutions found in Theorem \ref{T} is given by a kind of local
maximum for the functional related to problem $(P_\e)$, due to the 
concentrating coefficient $a$, while the other is a solution of saddle
type and is given by the interaction between the boundary of $\Omega$
and the concentration of the potential.
Since there is no assumption on the
topological complexity of $\Omega$ (for example $\Omega$ can be an
half-space),  in general it is not possible to say that problem 
$(P_\e)$ has more than 2 solutions.

This paper is organized as follows: in Section 2 we introduce a
suitable variational setting related to problem $(P_\e)$, we recall some
known results and define some notations; 
Section 3 is devoted to a compactness Lemma;
in Section 4 we prove
some preliminary results, used in the proof of Theorem \ref{T}, that
is contained in Section 5.


\sezione{The variational framework and useful tools}


In order to simplify the notations, in the following we will assume
$x_0=0$ and $a_0=1$. 
	Moreover if $u\in H^1_0(D)$, where $D$ is an open set in $\R^N$, 
we will denote with the same symbol $u$ its extension to $\R^N$,
obtained by setting $u=0$ outside $D$.

For $D\subset\R^N$, let us set $D_\e=\{x\in\R^N\ :\ \e x\in D\}$.
	A simple computation shows that if $v$ solves problem
$(P_\e)$, then $u(x)=v(\e x)$ solves
\[
(\tilde P_\e) \hspace{1cm}
\left\{
        \begin{array}{ll}
         -\Delta u+(1+\alpha(x))u=u^{p-1} & \mbox{ in }\Omega_\e\\
         u>0                            & \mbox{ in }\Omega_\e\\
         u=0                            & \mbox{ on }\partial\Omega_\e.
        \end{array}
\right.
\]

The solutions of $(\tilde P_\e)$ correspond to the nonnegative
functions that are critical points of the functional
$E:H_0^1(\Omega_\e)\to\R$ given by 
\[
E(u)=\int_{\R^N}(|\D u|^2+(1+\alpha(x))u^2)dx,
\]
constrained on the manifold
\[
V_\e=\{u\in H_0^1(\Omega_\e)\ :\ |u|_{L^p}=1\}.
\]

Let us set
\begin{equation}
\label{me}
m_\e=\inf\{E(u)\ :\ u\in V_\e\},
\end{equation}
\begin{equation}
\label{m0}
m=\min_{u\in H^1(\R^N),{|u|_{L^p}=1}}\int_{\R^N}(|\D u|^2+u^2)dx.
\end{equation}

	It is well known that the minimum in (\ref{m0}) is
achieved by a positive radial function $w$, that is decreasing when the radial
co-ordinate increases, unique modulo translation and such that
\begin{equation}
\label{as}
\lim_{|x|\to\infty}|D^j w(x)||x|^{\frac{N-1}{2}}e^{|x|}=d>0,
\hspace{3mm} d\in\R
\end{equation}
(see \cite{BeL} and \cite{GNN}).

\begin{prop}
\label{P1}
If $\Omega$ verifies $(C_1)$, then $m_\e= m$ for every $\e>0$ and the
minimization problem (\ref{me}) has no solution.
\end{prop}
 \noindent{\underline {\sf Proof}} \hspace{2mm}
Fix $\e>0$. Since $H_0^1(\Omega_\e)\subset H^1(\R^N)$, it follows at 
once $m_\e\geq m$.

To prove that, actually, equality holds, fix
$0<\nu<\dist(0,\partial\Omega)$ and set
\begin{equation}
\label{O-}
\Omega^-=\{x\in\Omega\ :\
\dist(x,\partial\Omega)>\nu\}
\end{equation}
Then define a smooth function $\zeta_\e$ on $\Omega_\e$  such that
\begin{equation}
\label{P13}
\supp\zeta_1\subset\Omega,\ \zeta_1(x)=1\mbox{ on }\Omega^-,\
\zeta_\e(x)=\zeta_1(\e x).
\end{equation}
Now, for every $n\in\N$ let us choose $x_n\in\Omega_\e$ such that
$B(x_n,n)\subset\Omega_\e^-$ (it is possible by $(C_1)$) and consider the functions
\[
v_n(x)=\frac{\zeta_\e(x)w(x-x_n)}{|\zeta_\e(x)w(x-x_n)|_{L^p}}.
\]
Taking into account (\ref{as}), it is not difficult to see that
\[
\lim_{n\to\infty}E(v_n)=m,
\]
so $m_\e=m$.

\vspace{1mm}

        Let us now assume that the minimization problem (\ref{me})
has a solution $u^*\geq 0$.
        Then
\[
m\leq \|u^*\|^2_{H^1(\R^N)}=\|u^*\|^2_{H^1(\Omega_\e)}\leq \|u^*\|^2_{H^1(\Omega_\e)}
+\int_{\Omega_\e}\alpha(x)(u^*(x))^2dx=m.
\]
Thus we deduce
\[
u^*(x)=w(x-y^*) \hspace{5mm} \mbox{ for some } y^*\in\R^N
\]
and, by $(\ref{a1})$ and $w(x)>0$ $\forall x\in\R^N$,
\[
0=\int_{\Omega_\e}\alpha(x)(u^*(x))^2dx=\int_{\Omega_\e}\alpha(x)w^2(x-y^*)dx
>0,
\]
that is a contradiction.

\hfill{\em q.e.d.}

\vspace{2mm}
   		
\vspace{2mm}

In the following of this section some tools are introduced, useful to
describe some topological properties of the sublevels of $E$ on $V_\e$.

First, a barycenter type function is defined.
For $u\in L^p(\R^N)$  set
\[
\tilde u(x)={1\over \omega_N}\int_{B(x,1)} |u(y)|dy,
\]
$\omega_N$ being the Lebesgue measure of the unit ball in $\R^N$, and 
\[
\hat u(x)=\left[ \tilde u(x)-{1\over 2}\max_{\R^N} \tilde u
(x)\right]^+;
\]
then define $\beta:L^p(\R^N)\setminus\{0\}\to\R^N$ by
\begin{equation}
\label{beta}
\beta(u)={1\over |\hat u|^p_{L^p}}\int_{\R^N} x\cdot[\hat u(x)]^pdx.
\end{equation}
Observe that $\beta$ is well defined for all $u\in
L^p(\R^N)\setminus\{0\}$, because $\hat u\not\equiv 0$ and has compact
support, moreover $\beta$ is continuous and verifies
\begin{equation}
\label{b0}
\beta(v(x))=\beta(v(x-z))-z\qquad\forall v\in L^p(\R^N)\setminus\{0\},\
\ \forall z\in\R^N.
\end{equation}

Then some functions in $V_\e$  and sets in $\R^N$ are defined.
Let us denote $w_\e:\R^N\to V_\e$ by
\[
w_\e[y](x)=\frac{\zeta_\e(x)w(x-y)}{|\zeta_\e(x)w(x-y)|_{L^p}},\quad
x,y\in\R^N,
\]
where $w$ is the minimizing function for $m$ (see (\ref{m0})) and
$\zeta_\e$ has been introduced in (\ref{P13}). 

Fixed $\xi$ such that $|\xi|=1$, set
\[
\Sigma:=\partial B(\xi,2)=\{y\in\R^N\ :\ |y-\xi|=2\},
\]
and, for $\e,\rho>0$, define $w_{\e,\rho}:\Sigma\times [0,1]\to V_\e$ by
\[
w_{\e,\rho}[y,t](x)=\frac{\zeta_\e(x)[(1-t)w(x-\rho y)+tw(x-\rho \xi)]}
                     {|\zeta_\e(x)[(1-t)w(x-\rho y)+tw(x-\rho \xi)]|_{L^p}}.
\]
In particular, it holds $w_{\e,\rho}[y,0]=w_\e[\rho y]$ and 
$w_{\e,\rho}[y,1]=w_\e[\rho \xi]$. 

Finally, let us consider $\bar x\in\partial \Omega$ such that $\{t\bar
x\ :\ t\in[0,1)\}\subset\Omega$ and call
\begin{equation}
\label{Se}
S_\e=\{t\bar x/\e\ :\ t\in[0,1-\e]\}.
\end{equation}


\sezione{A compactness results}


In this section we prove a result which  states
that, under conditions $(C_1)$ and $(C_2)$, the 
functional $E$ constrained on $V_\e$ verifies the Palais-Smale
condition in an energy range. This statement will allow us to apply
some mini-max techniques of the Calculus of Variations, in order to
prove Theorem \ref{T}.

\begin{lemma}
\label{PS}
Suppose that $\Omega$ verifies assumptions $(C_1)$ and $(C_2)$, fix $\e>0$ and let $(u_n)_n$ be 
a Palais-Smale sequence for $E$ constrained on
$V_\e$, i.e. $u_n\in V_\e$ and 
\begin{equation}
\label{PS1}
        {{\lim_{n\to\infty}{E}(u_n)=c }}
\end{equation}
\begin{equation}
\label{PS2}
\lim_{n\to\infty} \D {E}_{\mid{V_\e}}(u_n)= 0.
\end{equation}
If $c\in (m,2^{1-2/p}m)$ then $(u_n)_n$ is relatively compact in $V_\e$.
\end{lemma}

\proof To get Lemma \ref{PS}, the behaviour at infinity of the
(PS)-sequence $(u_n)_n$ is analysed, as in \cite{BC,L}. Roughly
speaking, it is proved that a sequence $(u_n)_n$ which verifies
(\ref{PS1}) and (\ref{PS2}) may be decomposed in ``waves'' at infinity
which solve some limit problems. Then, by the energy estimate
$c\in(m,2^{1-2/p}m)$, we can conclude that these waves must vanish,
hence $(u_n)_n$ converge strongly to its weak limit.

\vspace{1mm}

From (\ref{PS1}) it follows that $(u_n)_n$ is bounded in
$H_0^1(\Omega_\e)$, so there exists $v_0\in H^1_0(\Omega_\e)$ such
that, up to a subsequence, 
\begin{equation}
\label{PS2.5}
u_n\to v_0,\ \mbox{ as }n\to \infty,\ \mbox{ in } L^p_{\loc}(\Omega_\e)
\mbox{ and weakly in } H^1_0(\Omega_\e).
\end{equation}

If we consider $\bar u_n=u_n-v_0$, then we can apply Lemma 3.2 of
\cite{MP98} to $(\bar u_n)_n$ and we find a sequence $(v_k)_k$ in
$H^1(\R^N)$ and, for each $k\in\N$, a sequence $(y_n^k)_n$ in
$\Omega_\e$ such that
\begin{equation}
\label{PS3}
\lim_{n\to \infty}|y_n^k|=\infty\ \ \forall k\in\N,\
\lim_{n\to\infty}|y_n^{k'}-y_n^{k''}|=\infty\ \ \mbox{ if } k'\neq k''
\end{equation}
and
\begin{equation}
\label{PS4}
\bar u_n(x+y_n^k)\to v_k\ \mbox{ as }n\to\infty,\ \mbox{ in }L^p_{\loc}(\R^N)
\  \mbox{ and weakly in } H^1(\R^N);
\end{equation}
moreover
\begin{equation}
\label{PS5}
\sum_{k=1}^\infty |v_k|^p_{L^p}=1-|v_0|^p_{L^p},
\end{equation}
\begin{equation}
\label{PS6}
E(v_0)+\sum_{k=1}^\infty\|v_k\|_{H^1(\R^N)}\leq \lim_{n\to\infty}E(u_n)=c.
\end{equation}

Observe that (\ref{PS2}) implies that there exists a sequence
$(\mu_n)_n$ in $\R$ such that for every $v\in H_0^1(\Omega_\e)$
\begin{equation}
\label{PS7}
\int_{\Omega_\e}[\D u_n\D
v+(1+a(x))u_nv]dx=\mu_n\int_{\Omega_\e}|u_n|^{p-2}u_nvdx+o(1)\|v\|_{H^1(\Omega_\e)}.
\end{equation}
Therefore, putting $v=u_n$ in (\ref{PS7}) and using (\ref{PS1}), we
obtain 
\begin{equation}
\label{PS8}
\lim_{n\to\infty}\mu_n=c.
\end{equation}
In particular, from (\ref{PS2.5}), (\ref{PS7}) and (\ref{PS8}) it
follows that $v_0$ solves
\begin{equation}
\label{PS8.5}
\left\{
\begin{array}{ll}
-\Delta u+(1+a(x))u=c |u|^{p-2}u & \mbox{ in }\Omega_\e\\
u=0& \mbox{ on }\partial\Omega_\e.
\end{array}\right.
\end{equation}

Moreover, taking into account (\ref{PS4}), (\ref{PS7}),
(\ref{PS8}), $(C_1)$ and $(C_2)$, from standard arguments it follows that $v_k$ solves
\begin{equation}
\label{PS9}
\left\{
\begin{array}{ll}
-\Delta u+u=c |u|^{p-2}u  & \mbox{ in }\R^N\\
u\in H^1(\R^N)
\end{array}\right.
\end{equation}
or gives a solution of
\begin{equation}
\label{PS10}
\left\{
\begin{array}{ll}
-\Delta u+u=c |u|^{p-2}u  & \mbox{ in }\Pi\\
u\in H^1_0(\Pi),
\end{array}\right.
\end{equation}
where $\Pi$ is an half-space in $\R^N$. 

If $v_k$ gives a solution of (\ref{PS10}), then $v_k\equiv 0$, by Theorem I.1 in \cite{EL}.

\vspace{1mm}

It is well known that $u=(m/c)^{1/(p-2)} w$ (see (\ref{m0})) is the
least energy solution among all the nontrivial solutions of
(\ref{PS9}) (see \cite{BC}, for example), hence if $v_k$ solves
(\ref{PS9}) then
\begin{equation}
\label{PS11}
\|v_k\|_{H^1(\R^N)}\geq (m/c)^{2/(p-2)}m.
\end{equation}
Inequalities (\ref{PS11}) and (\ref{PS6}) imply
that there exists $\bar k\in\N$ such that $v_k= 0$ for $k\geq
\bar k$.
	Moreover if $v_{k'},v_{k''}\neq 0$, for $k'\neq k''$, then from
(\ref{PS11}) and (\ref{PS6}) it follows that $c\geq 2^{1-2/p}m$,
contrary to our assumption.
	Hence we must have at most $v_1\neq 0$, for example.

We claim that also $v_1=0$. Arguing by contradiction, assume that
$v_1\neq 0$. 

Observe that, by (\ref{PS8.5}), $c^{1/(p-2)}v_0$ solves $(\tilde
P_\e)$ (except for the sign), moreover all the
nontrivial solutions of $(P_\e)$ are of the form $[E(\bar
u)]^{1/(p-2)}\bar u$, where $\bar u$ is a critical point for $E$ on
$V_\e$. So, taking into account Proposition \ref{P1}, we get
\begin{equation}
\label{PS13}
E(v_0)>(m/c)^{2/(p-2)}m.
\end{equation}

From (\ref{PS13}), (\ref{PS11}) and (\ref{PS6}) it follows
that if $v_0,v_1\neq 0$ then $c>2^{1-2/p}m$, contrary to our
assumption. 
	Hence $v_1\neq 0$ implies $v_0=0$ and so by (\ref{PS4}) and
(\ref{PS5}) we get
$|v_1|_{L^p}=1$ and $u_n(x+y_n^1)\to v_1$, as $n\to\infty$, strongly
in $L^p$.
	Applying this result to (\ref{PS7}) and using (\ref{PS4}) we
obtain $c=\|v_1\|_{H^1(\R^N)}$. If $v_1(x)\geq 0$ $\forall x\in \R^N$
(or $v_1\leq 0$), then, as stated in \S 2, $\|v_1\|_{H^1(\R^N)}=m$,
while we have assumed $c>m$, so it has to be $v^+\not\equiv 0$ and
$v^-\not\equiv 0$.
	Then by (\ref{PS9}) and by the definition of $m$ (see (\ref{m0}))
\begin{equation}
\label{PS14}
m|v_1^\pm|^2_{L^p}\leq\int_{\R^N}(|\D
v_1^\pm|^2+(v_1^{\pm})^2)dx=c|v_1^\pm|^p_{L^p},
\end{equation}
hence $|v_1|_{L^p}^p\geq 2(m/c)^{p/(p-2)}$, that implies $c>2^{1-2/p}m,$
against our assumption.
	Therefore it must be $v_1=0$.

\vspace{1mm}

Finally, let us prove that $u_n\to u_0$ strongly in $H^1(\Omega_\e)$.
First observe that $u_n\to v_0$ strongly in $L^p(\Omega_\e)$, by
(\ref{PS2.5}) and (\ref{PS5}). Then let us compute
\[
\hspace{-3cm}\|u_n-v_0\|_{H^1(\Omega_\e)}<E(u_n-v_0)=E(u_n)+E(u_0)
\]
\begin{equation}
\hspace{2cm}-2\int_{\Omega_\e}[\D u_n\D
v_0+(1+(x))u_nv_0]dx\longrightarrow 0,\ \mbox{ as }n\to\infty,
\end{equation}
by (\ref{PS2.5}) and (\ref{PS7}), applied with $v=u_n$ and $v=v_0$.

\qed

\begin{rem}
\label{SL}
{\em
In the proof of Lemma \ref{PS} we have used assumptions $(C_1)$ and
$(C_2)$ in
order to obtain that the ``waves'' $v_k$, to which (PS)-sequences converge, 
give a solution of  one of the limit problems (\ref{PS8.5}), (\ref{PS9}) and (\ref{PS10}). 
	If we drop either assumption $(C_1)$  or $(C_2)$, the alternatives 
(\ref{PS8.5})-(\ref{PS10}) are not
the only ones possible and, in fact, Lemma \ref{PS} could be false.
	
To explain the reason why Lemma \ref{PS} does not hold, in general, 
if condition $(C_1)$ is not verified, consider, for example, an open set $\hat\Omega$
which at infinity looks like a strip-like domain. 
	For a suitably large strip-like domain $\Omega$ there
are positive solutions $\hat u$ to problem $-\Delta u+u=u^{p-1}$ in
$H_0^1(\Omega)$, whose energy $E(\hat u/|\hat u|_{L^p})$ is in
the range $(m,2^{1-2/p}m)$ (see \cite{E}).
	Therefore it is not difficult to construct not relatively
compact (PS)-sequences for the functional $E$ on $\{u\in
H_0^1(\hat\Omega)\ :\ |u|_{L^p}=1\}$ at the level $E(\hat u/|\hat
u|_{L^p})$.

Suppose now that assumption  $(C_2)$ is not fulfilled and let us show
that also in this case (PS)-condition fails, for some domains. 
	To this end, consider an open set $\tilde\Omega$ which at
infinity looks like the complementary of a cylinder.
	Also for the complementary $\Omega$ of a  cylinder
there are positive solutions $\tilde u$ to problem $-\Delta u+u=u^{p-1}$ in
$H_0^1(\Omega)$, whose energy $E(\tilde u/|\tilde u|_{L^p})$ is in
the range $(m,2^{1-2/p}m)$ (\cite{CPc}), hence we can conclude as
before.

Roughly speaking, (PS)-condition fails for domains like $\hat\Omega$
because they remain thin at infinity, while it fails for domains like
$\tilde\Omega$ because $\partial \tilde \Omega$  does not flatten at infinity
and it does not become smaller and smaller.
}
\end{rem}

The following lemma estabilishes a lower bound for the energy of a critical
point $u$ of $E$ on $V_\e$ which changes sign. 
The proof can be easily deduced using the definition of $m$, as in
(\ref{PS14}).

\begin{lemma}
\label{2.3}
Let $u\in H_0^1(\Omega_\e)$ be such that 
\[
|u|_{L^p}=1 \qquad E(u)=c \qquad \D {E}_{\mid{V_\e}}(u)=0.
\]
Then $u^+\not\equiv 0$ and $u^-\not\equiv 0$ implies $c>2^{1-2/p}m$.
\end{lemma}


\sezione{Analysis of some sublevels of $E$ on $V_\e$}


In the proof of Theorem \ref{T} the solution of maximum type related
to the concentrating coefficient $a$ has ``high energy'' 
and the saddle type solution has ``small energy''.
To find the saddle type solution solution, consider $r>0$ such that
$B(0,r)\subset\Omega^-$  (see (\ref{O-})) and define
\begin{equation}
\label{S}
\cS_{\e}=\inf\{E(u)\ :\ u\in V_\e,\ \beta(u)\in S_\e\},
\end{equation}
\begin{equation}
\label{S0}
\cS_{\e,0}=\max_{\partial B(0,r/\e)} E(w_\e[y]),
\end{equation}
\begin{equation}
\label{S1}
\cS_{\e,1}=\inf\{E(u)\ :\ u\in V_\e,\ \beta(u)\in\{0,(1-\e)\bar x/\e\}\}.
\end{equation}

To find the solution related to $a$, let us define
\begin{equation}
\label{A1}
\cA_{\e}=\inf\{E(u)\ :\ u\in V_\e,\ \beta(u)=0\},
\end{equation}
\begin{equation}
\label{A2}
\cA_{\e,\rho,1}=\max_{\Sigma\times [0,1]}E(w_{\e,\rho}[z,t]),
\end{equation}
\begin{equation}
\label{A3}
\cA_{\e,\rho,0}=\max_{\Sigma}E(w_{\e,\rho}[z,0]).
\end{equation}

Next results state some properties of the levels just defined.

\begin{lemma}
\label{L1}
If $D$ is a compact subset of $\Omega$, then
\begin{equation}
\label{L11}
\lim_{\e\to 0}\max_{y\in D}|\beta(w_\e[y/\e])-y/\e|=0.
\end{equation}
\end{lemma}

 \noindent{\underline {\sf Proof}} \hspace{2mm}
For $\e>0$, let $y_\e\in D$; we are proving that
\begin{equation}
\label{L12}
\lim_{\e\to0}|\beta(w_\e[y_\e/\e])-y_\e/\e|=0.
\end{equation}
Denote $w_\e^t[z](x)=w_\e[z+y_\e/\e](x+y_\e/\e)$, $x,z\in\R^N$; 
from (\ref{b0}) it follows that (\ref{L12}) is equivalent to 
\begin{equation}
\label{L14}
\lim_{\e\to 0}|\beta(w_\e^t[0])|=0.
\end{equation}

Taking into account (\ref{P13}), we get $0\leq w_\e^t[0](x)\leq 2w(x)$
for small $\e$; then, in particular,
\begin{equation}
\label{L15}
\lim_{\e\to 0}w_\e^t[0](x)=w(x)\quad\mbox{ in } L^p(\R^N).
\end{equation}
By the symmetry of $w$ we have $\beta(w(x-z))=z$, $\forall z\in \R^N$;
hence (\ref{L14}) follows from (\ref{L15}) and from the continuity of
$\beta$ with respect to the $L^p$-norm.

Now, a suitable choice of $y_\e$ yields (\ref{L11}).

\hfill{\em q.e.d.}

\vspace{2mm}
   		
\begin{rem}
\label{R1}
{\em
Lemma \ref{L1} implies that if $D\subset\Omega^-$ and if $U$ is a
neighborhood of $D$, then for small $\e$
the map
\[
z\mapsto\beta(w_\e[z])
\]
is homotopic in $U_\e$ to the identity map, by the
homotopy $\cK:[0,1]\times D_\e\to U_\e$ defined by
\begin{equation}
\label{R12}
\cK(\theta,z)=\theta\beta(w_\e[z])+ (1-\theta)z \qquad
0\leq\theta\leq 1.
\end{equation}
}
\end{rem}

\begin{prop}
\label{P2}
If $a_\e$ is of the form (\ref{a}) and verifies (\ref{a1}), then there
exists $\mu_{\alpha,1}>m$ such that
\begin{equation}
\label{P21}
\cA_{\e}>\mu_{\alpha,1}\qquad\forall\e>0.
\end{equation}
\end{prop}
 
\noindent{\underline {\sf Proof}} \hspace{2mm}
Arguing by contradiction, let us assume that there exist 
sequences $(\e_i)_i$ in $\R^+$, $u_i\in H^1_0(\Omega_{\e_i})$ such that
\begin{equation}
\label{P22}
|u_i|_{L^p}=1,\ \beta(u_i)=0
\end{equation}
\begin{equation}
\label{P23}
\lim_{i\to\infty}\int_{\Omega_{\e_i}}[|\D u_i|^2+(1+\alpha(x))u_i^2]dx=m.
\end{equation}
Moreover we can assume $u_i\geq 0$ in $\R^N$, $\forall i\in\N$.

Our first claim is that it must be $\e_i\geq c$, for a suitable
constant $c>0$.
	Assume, contrary to our claim, that $\e_i\to 0$, as
$i\to\infty$, up to a subsequence.
From (\ref{P23}) and (\ref{m0}) it follows that $u_i$ is a minimizing
sequence for $m$, hence, by the uniqueness of the solution of
(\ref{m0}), a sequence of points $(z_i)_i$ in $\R^N$ and a sequence of 
functions $\varphi_i$ in $H^1(\R^N)$ exist such that 
\begin{equation}
\label{P24}
u_i(x)=w(x-z_i)+\varphi_i(x) \qquad x\in \R^N
\end{equation}
with
\begin{equation}
\label{P25}
\lim_{i\to\infty}\varphi_i(x)=0 \qquad \mbox{ in } H^1(\R^N) \mbox{ and
in } L^p(\R^N).
\end{equation}
The same arguments of Lemma \ref{L1} show that
$\lim_{i\to\infty}(\beta(u_i)-z_i)=0$ hence, by (\ref{P22})
\begin{equation}
\label{P26}
\lim_{i\to\infty}z_i=0.
\end{equation}
Using (\ref{P22})--(\ref{P26}), we get
\begin{equation}
\label{P27}
m=\int_{\R^N}[|\D w|^2+w^2]\leq\int_{\R^N}[|\D
w|^2+(1+\alpha(x))w^2]=m,
\end{equation}
that implies $\int_{\R^N}\alpha w=0$. This is not possible since $w>0$
on $\R^N$ and $\alpha$ verifies (\ref{a1}), hence the claim follows.

\vspace{1mm}

If $\e_i\geq c$, $\forall i\in \N$, set $\tilde\Omega_c=\cup_{\e\geq
c}\Omega_\e$.
	By (\ref{P23}) and arguing as in Proposition \ref{P1}, we can
conclude that $(u_i)_i$ is a minimizing sequence for $E$ constrained
on $\{u\in H_0^1(\tilde\Omega_c)\ :\ |u|_{L^p}=1\}$.
	Then (\ref{P24}) and (\ref{P25}) hold for suitable sequences
$(z_i)_i$ in $\R^N$ and $(\varphi_i)_i$ in $H^1(\R^N)$.
	Moreover (\ref{P26}) is verified and so $\varphi_i\to- w$ in
$\R^N\setminus\tilde\Omega_c$. 
	This is in contradiction with (\ref{P25}) and shows that also
$\e_i\geq c$ is not possible.

\hfill{\em q.e.d.}

\vspace{2mm}
   		
\begin{prop}
\label{P3}
If $a$ is of the form (\ref{a}) and verifies (\ref{a1}),(\ref{a2}),
then there exist $\mu_{\alpha,2}\in (m,2^{1-2/p}m)$, $\rho_\alpha>0$
and $\e_1>0$ such that if $\e\in(0,\e_1)$ then
\begin{equation}
\label{P31}
\cA_{\e,\rho_\alpha,1}\leq\mu_{\alpha,2},
\end{equation}
\begin{equation}
\label{P32}
\cA_{\e,\rho_\alpha,0}<\mu_{\alpha,1};
\end{equation}
moreover
\begin{equation}
\label{P33}
\cA_{\e}\leq \cA_{\e,\rho_\alpha,1}.
\end{equation}
\end{prop}
\noindent{\underline {\sf Proof}} \hspace{2mm}
For $z\in\Sigma$ and $t\in[0,1]$, let us define
$w_{0,\rho}[z,t]:\R^N\to\R$ by
\[
w_{0,\rho}[z,t](x)=\frac{(1-t)w(x-\rho z)+tw(x-\rho \xi)}
                     {|(1-t)w(x-\rho z)+tw(x-\rho \xi)|_{L^p}}.
\]
As stated in Step 1 of Proposition 3.3 in \cite{CM}, there exist
$\rho_1>0$ and $\mu_{\alpha,2}\in(m,2^{1-2/p}m)$ such that if
$\rho>\rho_1$ then
\begin{equation}
\label{P35}
\max_{\Sigma\times[0,1]} 
\int_{\R^N} \left[|\D w_{0,\rho} [z,t] |^2+(1+\alpha
(x))(w_{0,\rho}[z,t])^2\right]dx<\mu_{\alpha,2}.
\end{equation}

Let us compute
\[
\int_{\R^N} \left[|\D w_{0,\rho} [z,0] |^2+(1+\alpha
(x))(w_{0,\rho}[z,0])^2\right]dx
\]
\[
=\int_{\R^N} \left[|\D w(x-\rho z)|^2+(1+\alpha
(x))w(x-\rho z)^2\right]dx
\]
\begin{equation}
\label{P35a}
=m+\int_{\R^N} \alpha(x)w(x-\rho z)^2 dx,
\end{equation}
with
\begin{equation}
\label{P35b}
\lim_{|y|\to\infty}\int_{\R^N} \alpha(x)w(x-y)^2 dx=0,
\end{equation}
by (\ref{a1}) and (\ref{as}).

Hence there exists $\rho_0$ such that 
\begin{equation}
\label{P36}
\max_{\Sigma} 
\int_{\R^N} \left[|\D w_{0,\rho} [z,0] |^2+(1+\alpha
(x))(w_{0,\rho} [z,0])^2\right]dx <\mu_{\alpha,1}\quad\forall\rho>\rho_0.
\end{equation}

 Let us fix
$\rho_\alpha>\max\{\rho_0,\rho_1\}$; we claim that   
for every compact set $K\subseteq\Sigma\times[0,1]$
\[
\hspace{-6cm}
\lim_{\e\to 0}\max_{(z,t)\in K} E\left(w_{\e,\rho_\alpha}[z,t]\right)=
\]
\begin{equation}
\label{P37}
\hspace{2cm}\max_{(z,t)\in K} \int_{\R^N}
(|\D w_{0,\rho_\alpha}[z,t]|^2+
(1+\alpha(x))(w_{0,\rho_\alpha}[z,t])^2)dx.
\end{equation}
To prove (\ref{P37}), let  $(\e_i)_i$ in $\R^+$ and $(z_i,t_i)_i$ in
$ K$ be sequences such that
$\e_i\to 0$ and $(z_i,t_i)\to(z_0,t_0)\in K$, as $i\to\infty$.
	Since
$\dist(\{\rho_\alpha\xi\}\cup\rho_\alpha\Sigma,\R^N\setminus\Omega_\e)\to\infty$
as $\e\to 0$ and taking into account (\ref{as}), it is not difficult
to see that
\begin{equation}
\label{P38}
\lim_{i\to \infty} w_{\e_i,\rho_\alpha}[z_i,t_i]=w_{0,\rho_\alpha}
[z_0,t_0] \qquad \mbox{ in } H^1(\R^N).
\end{equation}
Hence we get (\ref{P37}), which implies (\ref{P31}) and (\ref{P32}), in
view of (\ref{P35}) and (\ref{P36}).

To show (\ref{P33}), observe that by  Remark
\ref{R1} the homotopy
$\cK:[0,1]\times\Sigma\to\Omega_\e\setminus\{0\}$ given by
\begin{equation}
\label{P39}
\cK(\theta,z)=\theta \beta(w_{\e,\rho_\alpha}[z,0])+(1-\theta)
\rho_\alpha z
\end{equation}
is well defined, for small $\e$. 
	Then, as $\beta$ is continuous, there exists $(\hat z,\hat
t)\in\Sigma\times[0,1]$ such that $\beta(w_{\e,\rho_\alpha}[\hat z,\hat t])=0$, from which
(\ref{P33}) follows.

\hfill{\em q.e.d.}

\vspace{2mm}
   		
\begin{prop}
\label{M}
If $a_\e$ is of the form (\ref{a}) and verifies (\ref{a1}), then there
exist $\mu_{0}\in(m,\mu_{\alpha,1}]$ and $\e_2>0$ such that if
$\e\in(0,\e_2)$ then
\begin{equation}
\label{MM1}
\mu_0<\cS_{\e,1},
\end{equation} 
\begin{equation}
\label{MM2}
\cS_{\e,0}<\mu_0,
\end{equation} 
\begin{equation}
\label{MM3}
m<\cS_\e\leq\cS_{\e,0}.
\end{equation} 
\end{prop}
 \noindent{\underline {\sf Proof}} \hspace{2mm}
Since (\ref{P21}) holds, to obtain (\ref{MM1}) it is enough to prove
that if $(\e_i)_i$ is a vanishing sequence in $\R^+$ and $u_i\in V_{\e_i}$
verify
\begin{equation}
\label{MM4}
\lim_{i\to\infty} E(u_i)=m,
\end{equation}
then
\begin{equation}
\label{MM5}
\dist(\beta(u_i),\bar x/\e_i)=+\infty.
\end{equation}

We can assume $u_i\geq 0$ in $\R^N$, $\forall i\in\N$.
By the uniqueness of the minimizer for (\ref{m0}), there exist a
sequence of points $(y_i)_i$ in $\R^N$ and a sequence of functions $(\phi_i)_i$
in $H^1(\R^N)$ such that
\begin{equation}
\label{L35}
u_i(x)=w(x-y_i)+\varphi_i(x)
\end{equation}
with
\begin{equation}
\label{L36}
\lim_{i\to \infty}\varphi_i=0 \qquad\mbox{ in }H^1(\R^N)\ \mbox{ and in }
L^P(\R^N).
\end{equation}
As in Lemma \ref{L1}, we obtain
\begin{equation}
\label{MM6}
\lim_{i\to\infty}|\beta(u_i)- y_i|=0,
\end{equation}
hence (\ref{MM5}) is equivalent to 
\begin{equation}
\label{MM7}
\lim_{i\to\infty}|y_i-\bar x/\e_i|=+\infty.
\end{equation}
Observe that  (\ref{L35}) and (\ref{L36}) imply that
\begin{equation}
\label{L311}
\lim_{i\to\infty} w(x-y_i)=0\qquad\mbox{ in }
L^p(\R^N\setminus\Omega_{\e_i}),
\end{equation}
so (\ref{MM7}) must hold and we get (\ref{MM1}).

\vspace{1mm}

By (\ref{as}) and taking into account that $\partial
B(0,r/\e)\subset\Omega^-$  and (\ref{P13}), it is not
difficult to verify that
\begin{equation}
\label{P510}
\lim_{\e\to 0} w_\e[y](x)=w(x-y)\qquad\mbox{ in }H^1(\R^N)
\end{equation}
uniformly with respect to $y\in\partial B(0,r/\e)$. Then by (\ref{a1})
we have
\[
\lim_{\e\to 0}\cS_{\e,0}=m,
\]
and (\ref{MM2}) follows.

\vspace{1mm}

As in Remark \ref{R1}, we have that for small $\e$
\begin{equation}
\label{M10}
\cK:[0,1]\times\partial B(0,r/\e)\to\Omega_\e\setminus\{0,(1-\e)\bar
x/\e\}
\end{equation}
defines a continuous map.

Hence there exists $\bar y_\e\in\partial B(0,r/\e)$ such that
$\beta(w_\e[\bar y_\e])\in S_\e$, that proves the second inequality
in (\ref{MM3}).

To see that $\cS_\e>m$, it is enough to observe that, for every fixed
$\e$, if a sequence $(u_i)_i$ in $V_\e$ verifies $E(u_i)\to m$, as
$i\to\infty$, then there exist a sequence of points $(y_i)_i$ in
$\R^N$ and a sequence of functions $(\varphi_i)_i$ in
$H^1(\R^N)$ with $|y_i|\to\infty$ and $\varphi_i\to 0$ in
$H^1(\R^N)$ and in $L^p(\R^N)$, as $i\to \infty$, such
that
\[
|u_i(x)|=w(x-y_i)+\varphi_i(x).
\]
Then $|\beta(u_i)-y_i|\to 0$, as $i\to\infty$, and, for large $i$,
$\beta(u_i)\not\in S_\e$.

\hfill{\em q.e.d.}

\vspace{2mm}


\sezione{Proof of Theorem \ref{T}}


In this proof we consider $0<\e<\bar\e:=\min\{\e_1,\e_2\}$ 
and the constants $\mu_{0}$, $\mu_{\alpha,1}$ and
$\mu_{\alpha,2}$ previosly defined (see Propositions \ref{P2}, \ref{P3}  and \ref{M}).
Moreover, for
$c\in\R$, we will set
\[
E_\e^c=\{u\in V_\e\ :\ E(u)\leq c\}.
\]

\vspace{1mm}

\noindent {\underline {\em Step 1}} \hspace{3mm} {\em  
Solution of saddle type.}

\vspace{1mm}

We claim that there exists a critical level
$c_{0,\e}\in[\cS_\e,\cS_{\e,0}]$ for the function $E$ on $V_\e$.
If it is not true, by Proposition \ref{M} and Lemma \ref{PS} we can apply a
well known deformation Lemma (see f.i. \cite{St}) and find a number
$\delta_0$ and a continuous function $\eta_0:[0,1]\times
E_\e^{\cS_{\e,0}}\to E_\e^{\cS_{\e,0}}$ such that
\begin{eqnarray}
\label{TT1}
& \eta_0(0,u)=u & \quad\forall u\in E_\e^{\cS_{\e,0}},\\
\label{TT2}
& \eta_0(1,u)\in E_\e^{\cS_{\e}-\delta_0} & \quad\forall u\in E_\e^{\cS_{\e,0}}.
\end{eqnarray}
Then the deformation $\cG:[0,1]\times\partial
B(0,r/\e)\to\R^N\setminus\{0,(1-\e)\bar x/\e\}$ given by
\[
\cG(t,y)=\left\{\begin{array}{ll}
\cK(2t,y) & \mbox{ if } t\in [0,1/2]\\
\beta\circ\eta_0(2t-1,w_\e[y]) & \mbox{ if }t\in[1/2,1]
\end{array}
\right.
\]
(see (\ref{R12})) is well defined and continuous (see (\ref{M10}),
(\ref{MM1}) and (\ref{MM2})).

By (\ref{R12}), (\ref{TT1}) and (\ref{TT2}), the map $\cG$ provides a continuous
deformation in $\R^N\setminus \{0,(1-\e)\bar x/\e\}$ from $\partial
B(0,r/\e)$ into a set that does not intersect $S_\e$. This is not
possible, so we get the claim.
\pagebreak


\noindent {\underline {\em Step 2}} \hspace{3mm} {\em  
Solution related to the coefficient $a$.}

\vspace{1mm}

We claim that there exists a critical value $c_{1,\e}\in
[\cA_{\e},\cA_{\e,{\rho_\alpha},1}]$.
If this is not the case, by Propositions \ref{P2}, \ref{P3} and by Lemma
\ref{PS}, we can find a number $\delta_1>0$
and a continuous function 
$\eta_1:E_\e^{\cA_{\e,{\rho_\alpha},1}}\to
E_\e^{\cA_{\e}-\delta_1}$ such that
\begin{equation}
\label{T1}
 \eta_1 (u)=u \hspace{2cm} \forall u\in
E_\e^{\cA_{\e}-\delta_1},
\end{equation}
furthermore, by (\ref{P21}) and (\ref{P32}),  $\delta_1$ can be chosen
in such a way that
\begin{equation}
\label{T2}
\cA_{\e,\rho_\alpha,0}<\cA_{\e}-\delta_1.
\end{equation}

Setting 
\[
\widetilde \Sigma={\Sigma\times [0,1]\over \sim},
\]
where $\sim$ identifies the points $(z,1)$, we define  a map $\cJ$ on
$\widetilde \Sigma$ by
\[
\cJ[z,t]=\beta\circ\eta_1(w_{\e,\rho_\alpha}[z,t]).
\]

By (\ref{T1}), (\ref{T2}) and Proposition \ref{P3} (see (\ref{P39})), $\cJ$  maps 
$\partial\widetilde \Sigma$ in a set homotopically equivalent 
to  $\rho_1\Sigma$ (and then to $\Sigma$) in $\R^N\setminus\{0\}$.
	Moreover $\cJ$ is  continuous, so a point $(\tilde z,\tilde
t)\in \widetilde \Sigma$ must exist, for which
\[
0=\cJ(\tilde z, \tilde t)=\beta\circ\eta_1(w_{\e,\rho_\alpha}[\tilde z,\tilde t]).
\]
This is impossible since  $\cJ(\widetilde\Sigma)\subset
\beta\circ\eta_1(E_\e^{\cA_{\e,{\rho_\alpha},1}})=
\beta(E_\e^{\cA_{\e}-\delta_1  })$ and by the definition of
$\cA_{\e}$ (see (\ref{A1})), so we are in contradiction.

\vspace{2mm}

Finally, let us remark that the critical levels $c_{0,\e}$ and
$c_{1,\e}$ are distinct, because
\begin{equation}
\label{T6}
m<\cS_\e\leq c_{0,\e}\leq\cS_{\e,0}<\mu_{\alpha,1}<c_{1,\e}\leq\cA_{\e,\rho_\alpha,1}<2^{1-2/p}m
\end{equation}
by Propositions \ref{P2}, \ref{P3} and \ref{M},  hence we get two
distinct critical points for $E$ on $V_\e$. 

Furthermore, the solutions related to these critical points are positive by Lemma \ref{2.3}, by
(\ref{T6}) and by the maximum principle. 
This completes the proof.

\hfill{\em q.e.d.}



\end{document}